\documentclass[fleqn]{mat01}
\usepackage{times,mathtimy,amssymb,latexsym,amscd}
\begin{document}

\setcounter{page}{483} \firstpage{483}

\newtheorem{theore}{Theorem}
\renewcommand\thetheore{\arabic{theore}}
\newtheorem{theor}[theore]{\bf Theorem}
\newtheorem{lem}[theore]{\it Lemma}
\newtheorem{propo}[theore]{\rm PROPOSITION}
\newtheorem{coro}[theore]{\rm COROLLARY}
\newtheorem{definit}[theore]{\rm DEFINITION}
\newtheorem{probl}[theore]{\it Problem}
\newtheorem{exampl}[theore]{\it Example}

\title{$\pmb{D}$-boundedness and $\pmb{D}$-compactness in finite
dimensional probabilistic normed spaces}

\markboth{Reza Saadati and Massoud Amini}{Finite-dimensional
probabilistic normed spaces}

\author{REZA SAADATI$^1$ and MASSOUD AMINI$^2$}

\address{$^1$Department of Mathematics, Islamic Azad University~--~Ayatollah Amoli Branch, Iran and Institute for Studies in Applied Mathematics, 1, 4th Fajr, Amol~46176-54553, Iran\\
\noindent $^2$ Department of Mathematics, Tarbiat Modarres
University, P.O. Box 14115-175, Tehran, Iran\\
\noindent E-mail: rsaadati@eml.cc; mamini@modares.ac.ir}

\volume{115}

\mon{November}

\parts{4}

\pubyear{2005}

\Date{MS received 1 March 2005; revised 9 May 2005}

\begin{abstract}
In this paper, we prove that in a finite dimensional probabilistic
normed space, every two probabilistic norms are equivalent and we
study the notion of $D$-compactness and $D$-boundedness in
probabilistic normed spaces.
\end{abstract}


\maketitle

\section{Introduction and preliminaries}

K~Menger introduced the notion of a probabilistic metric space in
1942 and since then the theory of probabilistic metric spaces has
developed in many directions \cite{sch1}. The idea of Menger was
to use distribution functions instead of nonnegative real numbers
as values of the metric. The notion of a probabilistic metric
space corresponds to situations when we do not know exactly the
distance between two points, but we know only probabilities of
possible values of this distance. Such a probabilistic
generalization of metric spaces appears to be well adapted for the
investigation of physical quantities and physiological thresholds.
It is also of fundamental importance in probabilistic functional
analysis. Probabilistic normed spaces were introduced by
\v{S}erstnev \cite{ser} in 1962 by means of a definition that was
closely modelled on the theory of (classical) normed spaces, and
used to study the problem of best approximation in statistics. In
the sequel, we shall adopt the usual terminology, notation and
conventions of the theory of probabilistic normed spaces, as in
\cite{al1,al2,gu,gu1,gu3,kh,sch1,si}.

In the sequel, the space of probability distribution functions
(briefly, d.f.) is $ \Delta^{+} =\{ F\hbox{:}\ {\bf R}\rightarrow
[0,1]\hbox{:}\ F $ is left-continuous, nondecreasing,  $F(0)=0$
and $F(+\infty)=1\}$ and the subset $D^{+} \subseteq \Delta^{+}$
is the set
\begin{equation*}
D^{+}=\{F\in  \Delta^{+}\hbox{:}\ l^{-}F(+\infty)=1\}.
\end{equation*}
Here $l^- f(x)$ denotes the left limit of the function $f$ at the
point $x$, $l^- f(x)=\lim_{t\to x^-}f(t)$. The space $\Delta^{+}$
is partially ordered by the usual point-wise ordering of
functions, i.e., $F\leq G$ if and only if $F(x)\leq G(x)$, for all
$x$ in ${\bf R}$. The maximal element for $\Delta^{+}$ in this
order is the d.f. given by
\begin{equation*}
\varepsilon_{0} = \begin{cases}
0, &\hbox{if} \; x\leq 0,\\[.2pc]
1, &\hbox{if} \; x>0.
\end{cases}
\end{equation*}
A triangle function is a binary operation on $\Delta^{+}$, namely
a function $\tau\hbox{:}\ \Delta^{+}\times \Delta^{+}\rightarrow
\Delta^{+}$ that is associative, commutative, nondecreasing and
which has $\varepsilon_{0}$ as unit, that is
\begin{align*}
&\tau(\tau(F,G),H) = \tau(F,\tau(G,H)),\\[.2pc]
&\tau(F,G) = \tau(G,F),\\[.2pc]
&F\leq G \Longrightarrow \tau(F,H)\leq\tau(G,H),\\[.2pc]
&\tau(F,\varepsilon_{0}) = F
\end{align*}
for all $F,G,H\in\Delta^{+}$. Continuity of triangle functions
means continuity with respect to the topology of weak convergence
in $\Delta^{+}$.

Typical continuous triangle functions are $\tau_{T}(F,G)
(x)=\sup_{s+t=x}T(F(s),G(t))$ and $\tau_{T^{*}}(F,G)=
\inf_{s+t=x}T^{*}(F(s),G(t))$. Here $T$ is a continuous $t$-norm,
i.e., a continuous binary operation on $[0,1]$ that is
commutative, associative, nondecreasing in each variable and has 1
as identity; $T^{*}$ is a continuous $t$-conorm, namely a
continuous binary operation on $[0,1]$ which is related to the
continuous $t$-norm $T$ through $ T^{*}(x,y)=1-T(1-x,\break 1-y).$

The definition below is more general, it has been proposed in
\cite{al1}.

\begin{definit}$\left.\right.$\vspace{.5pc}

\noindent {\rm A \emph{probabilistic normed space} (briefly, a PN
space) is a quadruple $(V,\nu,\tau,\tau^*)$, where $V$ is a real
vector space, $\tau$ and $\tau^*$ are continuous triangle
functions with $\tau\leq\tau^*$ and $\nu$ is a mapping (the
\emph{probabilistic norm}) from $V$ into $\Delta^+$, such that for
every choice of $p$ and $q$ in $V$ the following hold:
\begin{enumerate}
\leftskip .7pc
\item[(N1)] $\nu_p=\varepsilon_0$  if and only if  $p=\theta$
($\theta$ is the null vector in $V$);

\item[(N2)] $\nu_{-p}=\nu_p$;

\item[(N3)] $\nu_{p+q}\geq\tau(\nu_p, \nu_q)$;

\item[(N4)] $\nu_p\leq\tau^*(\nu_{\lambda p},\nu_{(1-\lambda)p})$
for every $\lambda\in[0,1]$.
\end{enumerate}
A PN space is called a \emph{\v{S}erstnev space} if it satisfies
(N1), (N3) and the following condition:

For every $\alpha\neq 0 \in {\bf R}$ and $x>0$  one has
\begin{equation*}
\nu_{\alpha p}(x)=  \nu_p \left(\frac{x}{|\alpha|} \right),
\end{equation*}
which clearly implies (N2) and also (N4) in the strengthened form
\begin{equation*}
\forall \lambda\in [0,1],\quad \nu_p=\tau_M(\nu_{\lambda
p},\nu_{(1-\lambda)p}).
\end{equation*}
A PN space in which $\tau=\tau_T$ and $\tau^* =\tau_{T^*}$ for a
suitable continuous $t$-norm $T$ and its $t$-conorm $T^*$ is
called a \emph{Menger PN space}.}
\end{definit}

\begin{lem}\hskip -.3pc{\rm \cite{al2}.} \ \
If $|\alpha|\leq |\beta|${\rm ,} then $\nu_{\beta p}\leq
\nu_{\alpha p}$ for every $p$ in $V$.
\end{lem}

\begin{lem}
\hskip -.3pc{\rm \cite{al2}.} \ \ If $\tau^{*}$ is Archimedean{\rm
,} then for every $p$ in $V$ such that $\nu_{p}\neq
\varepsilon_{\infty}$ and every $h>0${\rm ,} there is a $\delta>0$
such that
\begin{equation*}
|\alpha|<\delta \Longrightarrow \nu_{\alpha p}(h)>1-h.
\end{equation*}
\end{lem}

In this case PN space is a topological vector space (shortly,
TVS). We call every PN space with the above properties as
\emph{strong TVS}.

\begin{definit}{\rm \cite{sch1}.}$\left.\right.$\vspace{.5pc}  

\noindent{\rm Let $(V,\nu,\tau,\tau^{*})$ be a PN space. For each
$p$ in $V$ and $\lambda>0$, the strong $\lambda$-{\it
neighborhood} of $p$ is the set
\begin{equation*}
N_{p}(\lambda) = \{q\in V\hbox{:}\ \nu_{p-q}(\lambda)>1-\lambda\},
\end{equation*}
and the strong neighborhood system for $V$ is the union
$\cup_{p\in V}\mathcal{N}_{p}$ where $\mathcal{N}_p=
\{N_p(\lambda)\hbox{:}\break \lambda>0\}$.}
\end{definit}

The strong neighborhood system for $V$ determines a Hausdorff
topology for $V$.

\begin{definit}{\rm \cite{sch1}.}$\left.\right.$\vspace{.5pc} 

\noindent {\rm Let $(V,\nu,\tau,\tau^{*})$ be a PN space, a
sequence $\{p_{n}\}_{n}$ in $V$ is said to be strongly convergent
to $p$ in $V$ if for each $\lambda>0$, there exists a positive
integer $N$ such that $p_{n}\in N_{p}(\lambda)$, for $n\geq N$.
Also the sequence $\{p_{n}\}_{n}$ in $V$ is called strongly Cauchy
sequence if for every $\lambda >0$ there is a positive integer $N$
such that $\nu_{p_{n}-p_{m}}(\lambda)>1-\lambda $, whenever
$m,n>N$. A PN space $(V,\nu,\tau,\tau^{*})$ is said to be strongly
complete in the strong topology if and only if every strongly
Cauchy sequence in $V$ is strongly convergent to a point in $V$.}
\end{definit}

\begin{definit}{\rm \cite{gu1}.}$\left.\right.$\vspace{.5pc} 

\noindent {\rm Let $(V,\nu,\tau,\tau^{*})$ be a PN space and $A$
be the nonempty subset of $V$. The probabilistic radius of $A$ is
the function $R_{A}$ defined on ${\bf R}^{+}$ by
\begin{equation*}
R_{A}(x) = \begin{cases}
l^{-}\inf\{\nu_{p}(x)\hbox{:}\ p\in A\}, &\text{if} \;\; x\in[0,+\infty),\\[.2pc]
1, &\text{if} \;\; x=+\infty.
\end{cases}
\end{equation*}}
\end{definit}

\begin{definit}{\rm \cite{gu1}.}$\left.\right.$\vspace{.5pc} 

\noindent {\rm A nonempty set $A$ in a PN space
$(V,\nu,\tau,\tau^{*})$ is said to be:
\begin{enumerate}
\renewcommand\labelenumi{(\alph{enumi})}
\leftskip .1pc
\item certainly bounded, if  $R_{A}(x_{0})=1$ for some $x_{0}\in
(0,+\infty)$;

\item perhaps bounded, if one has $R_{A}(x)<1$, for every $x\in
(0,+\infty)$ and  $l^{-}R_{A}(+\infty)=1$;

\item perhaps unbounded, if  $R_{A}(x_{0})>0$ for some $x_{0}\in
(0,+\infty)$ and  $l^{-}R_{A}(+\infty)\in(0,1)$;

\item certainly unbounded, if  $l^{-}R_{A}(+\infty)=0$, i.e., if
$R_{A}=\varepsilon_{\infty}$.
\end{enumerate}

Moreover, $A$ is said to be distributionally bounded, or simply
$D$-bounded if either (a) or (b) holds, i.e., $R_{A} \in D^{+}.$
If $R_{A}\in \Delta^{+} \!\setminus\!\!D^{+}$, $A$ is called
$D$-unbounded.}
\end{definit}

\begin{theor}[\cite{gu1}] 
A subset $A$ in the PN space $(V,\nu,\tau,\tau^{*})$ is
$D$-bounded if and only if there exists a d.f. $G\in D^{+}$ such
that $\nu_{p}\geq G$ for every $p\in A$.
\end{theor}

If $A\subset {\bf R}$ is $D$-bounded then in general $A$ is not
classically bounded.

\begin{exampl}
{\rm We consider the PN space $({\bf R},\nu,\tau,{\bf M})$, where
$\tau$ is a triangle function such that
$\tau(\varepsilon_{c},\varepsilon_{d})\leq \varepsilon_{c+d}$,
$\textbf{M}$ is the maximal triangle function and the
probabilistic norm $\nu\hbox{:}\ {\bf R}\rightarrow
\pmb{\Delta}^{+}$ is defined by $\nu_{p}=
\varepsilon_{\frac{|p|}{a+|p|}}$ for every $p$ in ${\bf R}$ and
for a fixed $a>0$, with $\nu_{p}(+\infty)=1$ (see Theorem~5 of
\cite{gu2}). With this norm, ${\bf R}$ is $D$-bounded because
$\nu_{p}\geq \varepsilon_{1}$.}
\end{exampl}

\begin{definit}$\left.\right.$\vspace{.5pc} 

\noindent {\rm We say that the probabilistic norm $\nu\hbox{:}\
{\bf R}\rightarrow \pmb{\Delta}^{+}$ has the Lafuerza Guill\'{e}n
property (briefly, the LG-property) if, for every $x>0$,
$\lim_{p\rightarrow \infty}\nu_{p}(x)=0$, or, equivalently,
$\lim_{p\rightarrow \infty}\nu_{p}=\varepsilon_{\infty}$.}
\end{definit}

\begin{exampl} 
{\rm The probabilistic norm in the last example does not have the
LG-property.}
\end{exampl}

\begin{exampl} 
{\rm The quadruple $({\bf R},\nu,\tau_{\pi},\tau^{*}_{\pi})$,
where $\nu\hbox{:}\ {\bf R} \rightarrow \pmb{\Delta}^{+}$ is
defined by
\begin{equation*}
\nu_{p}(x)=\begin{cases}
0, &\hbox{if} \;\; x=0, \\[.2pc]
\exp(-|p|^{1/2}), &\hbox{if} \;\; 0<x<+\infty, \\[.2pc]
1, &\hbox{if} \;\; x=+\infty,
\end{cases}
\end{equation*}
and $\nu_{0}=\varepsilon_{0}$ is a PN space (see \cite{al1}) but
is not \v{S}erstnev space and the probabilistic norm has the
LG-property.}
\end{exampl}

\begin{lem}
In a PN space $({\bf R},\nu,\tau,\tau^{*})$ in which the
probabilistic norm has the LG-property{\rm ,} if $A\subset {\bf
R}$ is $D$-bounded then it is classically bounded.
\end{lem}

\begin{proof} If $A\subset {\bf R}$ is $D$-bounded, there exists a
d.f. $G\in D^{+}$ such that $\nu_{p}\geq G,$  for every $p\in A$
but if $A$ is not classically bounded, then for every $k>0$ there
exists a $p\in A$ such that $|p|>k$. Hence $\lim_{p\rightarrow
\infty}\nu_{p}(x)=0$. Therefore for every $x\in(0,+\infty)$ we
have $G(x)=0$, which is a contradiction.\hfill $\Box$
\end{proof}

The converse of the above lemma is, in general, not true. See
Example~12. Here the only $D$-bounded set is the singleton
$\{\theta\}$.

\begin{theor}[\!] 
If the PN space $({\bf R},\nu,\tau,\tau^*)$ is a TVS then it is
complete.
\end{theor}

\begin{proof}
Let $\{p_{m}\}$ be strongly Cauchy sequence, for every $m,n\in N\;
m>n$ we have
\begin{equation*}
\lim_{m,n\rightarrow\infty}\nu_{p_{m}-p_{n}}=\varepsilon_{0}.
\end{equation*}
By TVS property we have
\begin{equation*}
\nu_{\lim_{m,n}(p_{m}-p_{n})}=\varepsilon_{0}=\nu_{0}.
\end{equation*}
Hence $\{p_{m}\}$ is a classical Cauchy sequence in ${\bf R}$;
therefore, it is convergent to $p\in {\bf R}$, i.e.,
$p_{m}-p\rightarrow 0$, since the PN space is TVS we have
$\lim_{m}\nu_{p_{m}-p}=\nu_{0}=\varepsilon_{0}.$\hfill $\Box$
\end{proof}

\section{Finite dimensional PN space}

In this section, we are interested in some properties of a finite
dimensional PN space, in particular we introduce the definition of
equivalent norms in a PN space.

\begin{theor}[\!] 
Let $\{p_{1},\dots,p_{n}\}$ be a linearly independent set of
vectors in a PN space $(V,\nu,\tau,\tau^{*})$ such that $\tau^{*}$
is Archimedean and $\nu_{p}\neq \varepsilon_{\infty}${\rm ,} for
every $p\in V$. Then there is a number $c\neq 0$ and there exists
a probabilistic norm $\nu'\hbox{\rm :}\ {\bf R}\rightarrow
\pmb{\Delta}^{+}$ on the real PN space $({\bf
R},\nu',\tau',\tau^{\prime *})$ where $\tau^{\prime *}$ is
Archimedean and $\nu'_{p}\neq \varepsilon_{\infty}${\rm ,} such
that for every choice of $n$ real scalars
$\alpha_{1},\dots,\alpha_{n}$ we have
\begin{equation}
\nu_{\alpha_{1}p_{1}+\cdots +\alpha_{n}p_{n}}\leq
\nu'_{c(|\alpha_{1}|+\cdots +|\alpha_{n}|)}.
\end{equation}
\end{theor}

\begin{proof} We write $s=|\alpha_{1}|+\cdots +|\alpha_{n}|$. If
$s=0$, all $\alpha_{j}$ are zero, so (1) holds. Let $s>0$. We
define $\mu_{p}=\nu_{s p}$ and $\mu'_{r}=\nu'_{s r}$. Then (1) is
equivalent to the following\break inequality,
\begin{align}
\mu_{\beta_{1}p_{1}+\cdots +\beta_{n}p_{n}}\leq \mu'_{c},\quad
&\beta_{j}=\alpha_{j}/s,\, \left( \sum_{j=1}^{n}|\beta_{j}| = 1
\right).
\end{align}
Hence it suffices to prove the existence of $c\neq 0$ and $\mu'$
such that (2) holds. Suppose otherwise, then there exists a
sequence $\{q_{m}\}$ of vectors
\begin{equation*}
q_{m}=\beta_{1}^{(m)}p_{1}+\cdots +\beta_{n}^{(m)}p_{n}, \ \left(
\sum_{j=1}^{n}|\beta_{j}^{(m)}|= 1 \right),
\end{equation*}
such that $\mu_{q_{m}}\rightarrow \varepsilon_{0}$ as
$m\rightarrow \infty$. Since $\sum_{j=1}^{n} |\beta_{j}^{(m)}|=1$,
we have $|\beta_{j}^{(m)}|\leq 1$. Hence, the sequence
$\{\beta_{1}^{(m)}\}$ has a convergent subsequence. Let
$\beta_{1}$ denote the limit of such a subsequence, and let
$\{q_{1,m}\}$ denote the corresponding subsequence of $\{q_{m}\}$.
By the same argument, $\{q_{1,m}\}$ has a subsequence
$\{q_{2,m}\}$ for which the corresponding sequence of real scalars
$\{\beta_{2}^{(m)}\}$ converges say to $\beta_{2}$. Continuing
this process, we obtain a subsequence $\{q_{n,m}\}$ of $\{q_{m}\}$
such that
\begin{equation*}
q_{n,m}=\sum_{j=1}^{n}\gamma_{j}^{(m)}p_{j}, \left(
\sum_{j=1}^{n}|\gamma_{j}^{(m)}|=1 \right)
\end{equation*}
and $\gamma_{j}^{(m)}\rightarrow \beta_{j}$ as $m\rightarrow
\infty$. Hence
\begin{equation*}
\lim_{m\rightarrow \infty}q_{n,m}=q:=
\sum_{j=1}^{n}\beta_{j}p_{j},
\end{equation*}
where $\sum_{j=1}^{n}|\beta_{j}|=1$, since
\begin{align*}
\mu_{q_{n,m}-q} &=
\mu_{\sum_{j=1}^{n}(\gamma_{j}^{(m)}-\beta_{j})p_{j}}\\[.2pc]
&\geq  \tau^{n}(\mu_{(\gamma_{1}^{(m)}-\beta_{1})p_{1}},\dots,
\mu_{(\gamma_{n}^{(m)}-\beta_{n})p_{n}})\rightarrow
\varepsilon_{0}
\end{align*}
as $m\rightarrow\infty$. Since $\{p_{1},\dots,p_{n}\}$ is linearly
independent and not all $\beta_{j}$'s are zero, we have $q\neq
\theta$. Since $\mu_{q_{m}}\rightarrow \varepsilon_{0}$, we have
$\mu_{q_{n,m}}\rightarrow \varepsilon_{0}$. But,
\begin{equation*}
\mu_{q}=\mu_{(q-q_{n,m})+q_{n,m}} \geq
\tau(\mu_{q-q_{n,m}},\mu_{q_{n,m}})\rightarrow\varepsilon_{0},
\end{equation*}
and hence $q=\theta$. This contradicts $q\neq \theta$.\hfill
$\Box$
\end{proof}

The following example shows that in the above theorem we need the
field ${\bf R}$ to be a strong TVS.

\begin{exampl}
{\rm Consider the PN space $({\bf R},\nu,\tau,\tau^{*})$ where
$\tau^{*}$ is Archimedean and $\nu_{p}\neq \varepsilon_{\infty}$.
By the above theorem there exists a $c\neq 0$ and a probabilistic
norm $\nu'\hbox{:}\ {\bf R}\rightarrow \pmb{\Delta}^{+}$ such that
$\nu_{p}\leq \nu'_{cp}$. If in the PN space $({\bf
R},\nu',\tau',\tau^{\prime *})$ $\lim_{m}\nu'_{p_{m}}<
\varepsilon_{0}$ whenever $p_{m}{\rightarrow}0$ in ${\bf R}$, then
for the sequence $\{2^{-n}\}$ we have $\nu_{2^{-n}}\leq
\nu'_{c2^{-n}}$ and consequently $\varepsilon_{0}<
\varepsilon_{0}$, which is a contradiction.}
\end{exampl}

From now on all the fields are strong TVS.

\begin{theor}[\!]
Every finite dimensional subspace $W$ of a PN space
$(V,\nu,\tau,\tau^{*})$ where $\tau^{*}$ is Archimedean and
$\nu_{p}\neq \varepsilon_{\infty}$ for every $p\in V${\rm ,} is
complete. In particular{\rm ,} every finite dimensional PN space
is complete.
\end{theor}

\begin{proof} Let $\{q_{m}\}$ be a strong Cauchy sequence in $W$.
Let $\dim W=n$ and $\{w_{1},\dots,w_{n}\}$ be a linearly
independent subset of $W$. Then each $q_{m}$ has a unique
representation of the form
\begin{equation*}
q_{m}=\alpha_{1}^{(m)}w_{1}+\cdots+\alpha_{n}^{(m)}w_{n}.
\end{equation*}
Since $\{q_{m}\}$ is a strong Cauchy sequence, for every $h>0$
there is a positive integer $N$ such that
\begin{equation*}
\nu_{q_{m}-q_{k}}(h)>1-h
\end{equation*}
whenever $m,k\geq N$. By the above theorem and Lemma~2, we have,
for every $j=1,2,\dots,n,$
\begin{align*}
1-h &< \nu_{q_{m}-q_{k}}(h)\\[.2pc]
&= \nu_{\sum_{j=1}^{n}(\alpha_{j}^{(m)}-\alpha_{j}^{(k)})w_{j}}(h)\\[.2pc]
&\leq \nu'_{c\sum_{j=1}^{n}|\alpha_{j}^{(m)}-\alpha_{j}^{(k)}|}(h)\\[.2pc]
&\leq \nu'_{c|\alpha_{j}^{(m)}-\alpha_{j}^{(k)}|}(h),
\end{align*}
where $c\neq 0$, $\nu'\hbox{:}\ {\bf R}\rightarrow
\pmb{\Delta}^{+}$ and $m,k\geq N$. This shows that each of the $n$
sequences $\{\alpha_{j}^{(m)}\}_{m}$ where $j=1,2,\dots,n$ is a
strong Cauchy in $({\bf R},\nu',\tau',\tau^{\prime *})$. Hence it
converges, say to $\alpha_{j}$. Now let us define $
q=\alpha_{1}w_{1}+\cdots+ \alpha_{n}w_{n}.$ Clearly, $q\in W$.
Furthermore,
\begin{align*}
\nu_{q_{m}-q} &= \nu_{\sum_{j=1}^{n}(\alpha_{j}^{(m)}-\alpha_{j})w_{j}}\\[.2pc]
&\geq \tau^{n}(\nu_{(\alpha_{1}^{(m)}-\alpha_{1})w_{1}},
\dots,\nu_{(\alpha_{n}^{(m)}-\alpha_{n})w_{n}}),
\end{align*}
tends to $\varepsilon_{0}$ whenever $m\rightarrow\infty$. That is
$q_{m}$ strongly converges to $q$. Hence $W$ is\break
complete.\hfill $\Box$
\end{proof}

\begin{definit}$\left.\right.$\vspace{.5pc} 

\noindent {\rm A probabilistic norm $\nu\hbox{:}\ V\rightarrow
\Delta^{+}$ is said to be equivalent to a probabilistic norm
$\mu\hbox{:}\ V\rightarrow \Delta^{+}$, if
$p_{m}\stackrel{\nu}{\longrightarrow} p$ is equivalent to
$p_{m}\stackrel{\mu}{\longrightarrow}p$.}
\end{definit}

In the following example two equivalent norms in probabilistic
normed space are given.

\begin{exampl}
{\rm We consider two PN space $(V,\nu,\tau,\tau^{*})$ with
probabilistic norm $\nu_{p}=\varepsilon_{\|p\|}$ where
$\tau(\varepsilon_{c}, \varepsilon_{d})= \varepsilon_{c+d},
(c>0,d>0)$ and $\tau\leq \tau^{*}$ (see Example~1 of \cite{gu1})
and $(V,\mu,\tau',{\bf M})$ with probabilistic norm
$\mu_{p}=\varepsilon_{\frac{\|p\|}{a+\|p\|}}, a>0$ (see Theorem~5
of \cite{gu2}). It is easy to see that the two probabilistic norms
are equivalent and so identity map $J\hbox{:}\ V\rightarrow V$ is
continuous.}
\end{exampl}

\begin{theor}[\!] 
If $V$ is a finite dimensional vector spaces{\rm ,} then every two
probabilistic norms $\nu$ of $(V,\nu,\tau,\tau^{*})$ and $\mu$ of
$(V,\mu,\tau',\tau^{\prime *})$ are equivalent{\rm ,} whenever
$\tau^{*}$ and $\tau^{\prime *}$ are Archimedean{\rm ,}
$\nu_{p}\neq \varepsilon_{\infty}${\rm ,} and $\mu_{p}\neq
\varepsilon_{\infty}${\rm ,} for every $p$ in $V$.
\end{theor}

\begin{proof} Let $\{v_{1},\dots,v_{n}\}$ be a linearly
independent subset of $V$. Let $p_{m}
\stackrel{\nu}{\longrightarrow} p$. We know that both $p_{m}$ and
$p$ have a unique representation as
\begin{equation*}
p_{m}=\alpha_{1}^{(m)}v_{1}+\cdots +\alpha_{n}^{(m)}v_{n},
\end{equation*}
and $p=\alpha_{1}v_{1}+\cdots +\alpha_{n}v_{n}.$ By Theorem~15 and
Lemma~2, we have
\begin{align*}
\nu_{p_{m}-p}=\nu_{\sum_{j=1}^{n}(\alpha_{j}^{(m)}-\alpha_{j})v_{j}}
\leq \nu'_{c\sum_{j=1}^{n}|\alpha_{j}^{(m)}-\alpha_{j}|} \leq
\nu'_{c|\alpha_{j}^{(m)}-\alpha_{j}|},
\end{align*}
where $c\neq 0$ and $\nu'\hbox{:}\ {\bf R}\rightarrow
\pmb{\Delta}^{+}$. Therefore $\lim_{m}\nu'_{c(|\alpha_{j}^{(m)}
-\alpha_{j}|)}= \nu'_{c(\lim_{m}|\alpha_{j}^{(m)}-\alpha_{j}|)}=
\varepsilon_{0}$, that is $\alpha_{j}^{(m)}{\rightarrow}
\alpha_{j}$ in ${\bf R}$. But
\begin{align*}
\mu_{p_{m}-p} = \mu_{\sum_{j=1}^{n}(\alpha_{j}^{(m)}
-\alpha_{j})v_{j}} \geq \tau^{\prime n} (\mu_{(\alpha_{1}^{(m)}-
\alpha_{1})v_{1}}, \dots,\mu_{(\alpha_{n}^{(m)}-
\alpha_{n})v_{n}}),
\end{align*}
so by continuity of $\tau'$ we have
$p_{m}\stackrel{\mu}{\longrightarrow}p$. By the same argument
$p_{m}\stackrel{\mu}{\longrightarrow}p$ implies
$p_{m}\stackrel{\nu}{\longrightarrow}p$.\hfill $\Box$
\end{proof}

In the next example we show that there are two PN spaces which are
not equivalent even in a finite dimensional probabilistic normed
space. Indeed, since $\tau_{M}$ is not Archimedean, the first PN
space in the next example is not strong TVS.

\begin{exampl}
{\rm We consider PN space $({\bf R},\nu,\tau_{\bf W},\tau_{\bf
M})$ where $\nu_{0}=\varepsilon_{0}$ and $\nu_{p}=\frac{1}{|p|+2}
\varepsilon_{0}+ \frac{|p|+1}{|p|+2}\varepsilon_{\infty}$ for
$p\neq 0$, we know $\tau_{M}$ is not Archimedean, (see \cite{al2})
and PN space $({\bf R},\nu,\tau,\tau^{*})$ with probabilistic norm
$\nu_{p}=\varepsilon_{|p|}$ where $\tau(\varepsilon_{c},
\varepsilon_{d})=\varepsilon_{c+d},(c>0,d>0)$ and $\tau\leq
\tau^{*}$ (see Example~1 of \cite{gu1}). Now, the sequence
$\{1/n\}$ in the first PN space is not convergent but in the
second it is convergent. Therefore the above PN spaces are not
equivalent.}
\end{exampl}

\section{$\pmb{D}$-bounded and $\pmb{D}$-compact sets in PN
spaces}

\begin{theor}[\!] 
Let $(V,\nu,\tau,\tau^{*})$ be a PN space in which
$\nu(V)\subseteq D^{+}$ and $D^+$ is invariant under $\tau${\rm ,}
i.e.{\rm ,} $\tau(D^{+}\times D^{+})\subseteq D^{+}$. If
$p_{m}\rightarrow p$ in $V$ and $A=\{p_{m}\hbox{\rm :}\ m\in{\bf
N}\}${\rm ,} then $A$ is a $D$-bounded subset of $V$.
\end{theor}

\begin{proof} Let $p_{m}\rightarrow p$. Then there exists a
positive integer $N$ such that for every $m\geq N$ we have
$\nu_{p_{m}-p}\geq G$, for each $G \in D^{+}$. Therefore
\begin{align*}
\nu_{p_{m}}&\geq \tau(\nu_{p_{m}-p},\nu_{p}) \geq
\tau(G,\nu_{p}).
\end{align*}
If we put $H=\min\{\nu_{p_{1}},\dots, \nu_{p_{N-1}},
\tau(G,\nu_{p})\}$, then $H\in D^{+}$ and $\nu_{p_{m}}\geq H$, for
every $m\in {\bf N}$. Hence $A$ is a $D$-bounded set.\hfill $\Box$
\end{proof}

Note that, in the Example~12 in which $\nu(V)\subseteq
\Delta^{+}\backslash D^+$, the sequence $\big\{\frac{1}{m}\big\}$
is convergent but $A=\big\{\frac{1}{m}\hbox{:}\ m\in {\bf
N}\big\}$ is not a $D$-bounded set.

\begin{definit}$\left.\right.$\vspace{.5pc} 

\noindent{\rm The PN space $(V,\nu,\tau,\tau^{*})$ is said to be
\emph{distributionally compact} (simply $D$-compact) if every
sequence $\{p_{m}\}_{m}$ in $V$ has a convergent subsequence
$\{p_{m_{k}}\}$. A subset $A$ of a PN space
$(V,\nu,\tau,\tau^{*})$ is said to be $D$-compact if every
sequence $\{p_{m}\}$ in $A$ has a subsequence $\{p_{m_{k}}\}$
convergent to a vector $p\in A$.}
\end{definit}

\begin{lem} 
A $D$-compact subset of a PN space $(V,\nu,\tau,\tau^{*})$ in
which $\nu(V)\subseteq D^{+}$ and $D^{+}$ is invariant under
$\tau${\rm ,} is $D$-bounded and closed.
\end{lem}

\begin{proof} Suppose that $A\subseteq V$ is $D$-compact. If $A$
is $D$-unbounded, it contains a $D$-unbounded sequence $\{p_{m}\}$
such that $\nu_{p_m}<\varepsilon_m$. This sequence could not have
a convergent subsequence, since a convergent sequence must be
$D$-bounded by Theorem~22. The closedness of $A$ is trivial.\hfill
$\Box$
\end{proof}

As in the classical case, a $D$-bounded and closed subset of a
(finite dimensional) PN space is not $D$-compact in general, as
one can see from the next examples.

\begin{exampl}
{\rm We consider quadruple $({\bf Q},\nu,\tau_{\pi}, \tau_{\bf
t_{\bf 2}})$, where $\pi(a,b)=a.b$, $t_{2}(a,b)=\frac{1}{1+
[((1/a)-1)^{2}+ ((1/b)-1)^{2}]^{1/2}}$, for every $a,b\in(0,1)$
and probabilistic norm $\nu_{p}(t)=\frac{t}{t+|p|^{1/2}}$. It is
straightforward to check that $({\bf Q},\nu,\tau_{\pi}, \tau_{\bf
t_{\bf 2}})$ is a PN space. In this space, convergence of a
sequence is equivalent to its convergence in ${\bf R}$. We
consider the subset $A=[a,b]\cap {\bf Q}$, where $a,b\in {\bf R}
\backslash{\bf Q}$. Since $R_{A}(t)=\frac{t}{t+(\max\{|a|,
|b|\})^{1/2}}$, then $A$ is a $D$-bounded set and since $A$ is
closed in ${\bf Q}$ classically, and so is closed in $({\bf
Q},\nu,\tau_{\pi},\tau_{\bf t_{\bf 2}})$. We know $A$ is not
classically compact in ${\bf Q}$, i.e., there exists a sequence in
${\bf Q}$ with no convergent subsequence in a classical sense and
so in $({\bf Q},\nu,\tau_{\pi},\tau_{\bf t_{\bf 2}})$. Hence $A$
is not $D$-\break compact.}
\end{exampl}

\begin{exampl}
{\rm We consider the PN space introduced in Example~9. With this
probabilistic norm, ${\bf R}$ is $D$-bounded and closed. But ${\bf
R}$ is not $D$-compact, because the sequence $\{2^{m}\}$ in ${\bf
R}$ does not have any convergent subsequence in this space.}
\end{exampl}

\begin{exampl}
{\rm The quadruple $({\bf R},\nu,\tau,\textbf{M})$, where
$\textbf{M}$ is the maximal triangle function, and the
probabilistic norm is a map $\nu\hbox{:}\ {\bf R}\rightarrow
\pmb{\Delta}^{+}$ such that $\nu_{0}=\varepsilon_{0}$,
$\nu_{p}=\varepsilon_{\frac{a+|p|}{a}}$ if $p\neq 0$, $(a>0)$, and
$\tau(\varepsilon_{c},\varepsilon_{d})\leq
\varepsilon_{c+d},(c>0,d>0)$, is a PN space (see Theorem~4 of
\cite{gu2}). If $A$ is a nonempty, classically bounded set in
${\bf R}$, then there exists $s>0$ such that for every $p\in A$ we
have $|p|\leq s$. Since $\nu_{p}\geq \varepsilon_{\frac{a+s}{a}}$,
$A$ is $D$-bounded. Also it is trivial that $A$ is closed. Now we
show that $A$ is not $D$-compact. Assume, if possible, that $A$ is
$D$-compact and $\{p_{m}\}$ be an arbitrary sequence in $A$ which
has a subsequence $\{p_{m_{k}}\}$ convergent to some $p$ in $A$,
then we have
\begin{equation*}
\lim_{k}\nu_{p_{m_{k}}-p}=\lim_{k}\varepsilon_{\frac{a+|p_{m_{k}}-p|}{a}}
\neq\varepsilon_{0}.
\end{equation*}
This implies that $p$ is not in $A$, a contradiction.}
\end{exampl}

\begin{theor}[\!]
Consider a finite dimensional PN space $(V,\nu,\tau,\tau^{*})${\rm
,} where $\tau^{*}$ is Archimedean{\rm ,} $\nu_{p}\neq
\varepsilon_{\infty}${\rm ,} and $\nu(V)\subseteq D^{+}$ and
$D^{+}$ is invariant under $\tau${\rm ,} for every $p\in V$ on the
real field $({\bf R},\nu',\tau',\tau^{\prime *})${\rm ,} where
$\nu'$ has the LG-property. Every subset $A$ of $V$ is $D$-compact
if and only if $A$ is $D$-bounded and closed.
\end{theor}

\begin{proof} By Lemma~24, $D$-compact subsets $A$ of $V$ are
$D$-bounded and closed, so we need to prove the converse. Let
$\dim V=n$ and $\{w_{1},\dots,w_{n}\}$ be a linearly independent
subset of $V$. We consider any sequence $\{q_{m}\}$ in $A$. Each
$q_{m}$ has a representation $ q_{m}=\alpha_{1}^{(m)}w_{1}+\cdots
+\alpha_{n}^{(m)}w_{n}$. Since $A$ is $D$-bounded so is
$\{q_{m}\}$, and so there exists a d.f. $G\in D^{+}$ such that
\begin{align*}
G &\leq \nu_{q_{m}}\\[.5pc]
&\leq \nu_{\alpha_{1}^{(m)}w_{1}+\cdots +\alpha_{n}^{(m)}w_{n}}.
\end{align*}
By Theorem~15 and Lemma~2, we have
\begin{align*}
&\leq \nu'_{c(|\alpha_{1}^{(m)}|+\cdots+
|\alpha_{n}^{(m)}|)}\\[.5pc]
&\leq  \nu'_{c|\alpha_{j}^{(m)}|}.
\end{align*}
Hence for each fixed $j$, the sequence $\{\alpha_{j}^{(m)}\}$ is
$D$-bounded and since $\nu'$ has the LG-property then by Lemma~13,
it is also classically bounded. Therefore for every $1\leq j\leq
n$ the sequence $\{\alpha_{j}^{(m)}\}$ has a convergent
subsequence converging to some $\alpha_{j}, j=1,\dots,n$. As in
the proof of Theorem~15, we can construct a subsequence
$\{r_{m}\}$ of $\{q_{m}\}$ which converges to
$r:=\sum_{j=1}^{n}\alpha_{j}w_{j}$. Since $A$ is closed, $r\in A$.
This shows that each sequence $\{q_{m}\}$ in $A$ has a convergent
subsequence in $A$. Hence $A$ is $D$-compact.\hfill $\Box$
\end{proof}

\section*{Acknowledgements}

The authors would like to thank the referee for giving useful
comments and suggestions for the improvement of this paper. It is
a pleasure to thank Prof.~B~Lafuerza Guill\'{e}n for helpful
suggestions about the subject of this paper. Also we would like to
thank Prof.~B~Schweizer for many helpful comments.

\end{document}